\documentclass{article}

\usepackage{latexsym,bm}
\usepackage{amsmath}  %如果要使用多充积分号，需要调用这个宏包
\usepackage{amssymb}
\usepackage{mathrsfs,ntheorem}

\newtheorem{theorem}{Theorem}[section]
\newtheorem{definition}[theorem]{Definition}
\newtheorem{lemma}[theorem]{Lemma}

\numberwithin{equation}{section}				%公式以节为编号

\DeclareMathOperator{\spans}{span}

\begin{document}

\title{A Correction on a Proof of a Combinatorial Property of the Set of Minimal Vectors in Root Lattices $\mathbb{A}_n$}

\author{Kong Xiaoran, Peking University.}

\date{\today}

\maketitle

\begin{abstract}
In this papar, we point out some mistakes in a proof of an important combinatorial property of $S(\mathbb{A}_n)$, the set of all minimal vectors of lattice $\mathbb{A}_n$, and correct them in the last section. This property plays an essential role in classifying perfect lattices in euclidean space.
\end{abstract}

\section{Introduction.}

\quad \, As a member of the famous family of root lattices, lattice $\mathbb{A}_n$ has some interesting properties. In the traditional way to describe $\mathbb{A}_n$, we usually embed $\mathbb{A}_n$ to the $(n+1)$-dimensional Euclidean space $\mathbb{R}^{n+1}$. Choose $\varepsilon_0, \varepsilon_1, \varepsilon_2, \cdots, \varepsilon_n$ as the orthonormal basis of $\mathbb{R}^{n+1}$, and let $H \subset \mathbb{R}^{n+1}$ be the subspace orthogonal to the vector $e = \varepsilon_0 + \varepsilon_1 + \cdots + \varepsilon_n$. We can define
$$\mathbb{A}_n = \mathbb{Z}^{n+1} \cap H = \left\lbrace\left. \sum_{k=0}^{n} a_k \varepsilon_k \right| \sum_{k=0}^{n} a_k = 0, a_i \in \mathbb{Z} \right\rbrace , $$
which is a full-rank lattice in the $n$-dimensional Euclidean space $H$. We can easily find out that $\{ e_i = \varepsilon_0 - \varepsilon_i \mid i=1,2,\cdots, n \}$ form a basis of $\mathbb{A}_n$, which is usually called Korkine-Zolotareff basis. 

Since all components of a vector in $\mathbb{A}_n$ add up to 0, there does not exist a vector of norm 1 in $\mathbb{A}_n$. Noticing that all $e_i$ are of norm\footnote{In this paper, the norm of a verctor means the square of its euclidean length.} 2, the shortest vectors in $\mathbb{A}_n$ should have norm 2. Hence, the set $S(\mathbb{A}_n)$, consisting of all minimal vector of $\mathbb{A}_n$, should be $\{\varepsilon_i - \varepsilon_j \mid 0 \leqslant i \neq j \leqslant n \}$, which could also be expressed as
$$ S(\mathbb{A}) = \{\pm e_1, \pm e_2, \cdots ,\pm e_n \} \cup \{ \pm(e_i-e_j) \mid 1 \leqslant i < j \leqslant n \}.$$

The set $S(\mathbb{A}_n)$ has a fantastic combinatorial property that every collection of $n$ independent vectors in $S(\mathbb{A}_n)$ spans the whole lattice $\mathbb{A}_n$, which can be shown by reduction to absurdity (See \cite{MatPerf}, Proposition 6.1.1). Surprisingly, the converse is also true! To classify perfect lattices with maximal index 1, J. Martinet recorded the interesting combinatorial property of $S(\mathbb{A}_n)$ in his famous textbook on perfect lattices in Euclidean Space \cite{MatPerf}:

\begin{theorem}[\cite{MatPerf}, Lemma 6.1.3] \label{Thm AnCombi}
Let $M$ be a free $\mathbb{Z}$-module of rank $n$ and let $\mathcal{S}$ be a family of elements of $M$ which satisfy the following five properties:  \par
1). $0 \notin \mathcal{S}$; \par
2). $x\in \mathcal{S} \Longrightarrow -x \in \mathcal{S}$; \par 
3). $\mathcal{S}$ is of rank $n$; \par 
4). $|\mathcal{S}| \geqslant n(n+1)$; \par
5). $n$ arbitrary independent elements of $\mathcal{S}$ generate $M$. \newline
Then $\mathcal{S}$ possesses exactly $\frac{n(n+1)}{2}$ pairs $\pm x$, from which one can extract a basis $(e_1, \cdots , e_n)$ for $M$ such that
$$\mathcal{S} = \{\pm e_i \mid 1 \leqslant i \leqslant n \} \cup \{ \pm(e_i-e_j) \mid 1 \leqslant i < j \leqslant n \}.$$
\end{theorem}

He also gave a proof on this proposition in Section 6.1 of \cite{MatPerf}. Historically, this proposition perhaps appeared first in the earlier paper \cite{K-Z An} of A. Korkine and G. Zolotareff, which is written in French and does not have English version. They used this proposition to classify perfect quadratic forms, which is equivalent to perfect lattices geometrically. However, both in \cite{K-Z An} and \cite{MatPerf}, the proofs of Theorem \ref{Thm AnCombi} include some mistakes. 

In this paper, we document the proof from \cite{MatPerf} in Section 2, and illustrate the mistakes of the proof in Section 3. In Section 4, we give our corrected proof of Theorem \ref{Thm AnCombi}.

\section{Original Proof.}

\quad \, In this section, we document the original proof of Theorem \ref{Thm AnCombi} in Chapter 6 of J. Martinet's book \cite{MatPerf}. We merely copy the proof literally here, with some comments we make as footnotes. We do not correct any mistakes in the proof while copying, which we postpone to the section \ref{secCorPf}. Both in the original textbook \cite{MatPerf} and our paper here, the components of an element with respect to a given basis are written as a column-vector rather than a row-vector. \\

In \cite{MatPerf}, the proof of Theorem \ref{Thm AnCombi} goes as: \\

Before entering into the proof of Theorem \ref{Thm AnCombi}, we introduce the notion
of a characteristic determinant: 

\begin{definition}[\cite{MatPerf}, Definition 6.1.4.] \label{Def 6.1.4}
Let $M$ be a free module of rank $n$ endowed with a basis $\mathcal{B} = (e_l, \cdots ,e_n )$, let $\mathcal{S} \supset \mathcal{B}$ be a finite subset of $M$, and let $s = |\mathcal{S}|$. A characteristic determinant for the pair $(\mathcal{B}, \mathcal{S})$ is a determinant of order $ r \leqslant n$ extracted from the $n \times s$ matrix of the components in $\mathcal{B}$ of the elements of $\mathcal{S}$.
\end{definition}

\begin{lemma}[\cite{MatPerf}, Lemma 6.1.5.] \label{Lemma 6.1.5.}
With the notation of Definition \ref{Def 6.1.4}, let $K$ be the maximum of the absolute values of determinants $\det_{\mathcal{B}}(x_1, \cdots ,x_n)$ for $x_1$, $\cdots$, $x_n \in \mathcal{S}$. Then the inequality $|d| \leqslant K$ holds for any characteristic determinant $d$ of $(\mathcal{B}, \mathcal{S})$. In particular, the absolute values of the components in $\mathcal{B}$ of any element of $\mathcal{S}$ are bounded from above by $K$.
\end{lemma}

\textit{Proof.} \; A characteristic determinant $d$ is the determinant of the matrix $D$ of $r \leqslant n$ elements $x_1, \cdots ,x_r$ of $S$ relatively to a system of $r$ vectors $b_1, \cdots , b_r$ of $\mathcal{B}$. After performing a suitable permutation of the $e_i$, we may assume that $b_i = e_i$. We then have
$$ \pm d = \operatorname{det}_{(e_1, \cdots ,e_r)} (x_1, \cdots ,x_r) = \operatorname{det}_{(e_1, \cdots ,e_n)} (x_1, \cdots ,x_r, e_{r+1}, \cdots ,e_n),$$
and the absolute value of the last determinant is bounded above by $K$. \hfill $\square$ \\

\textit{Proof of \ref{Thm AnCombi}.} We use induction on $n$, Theorem \ref{Thm AnCombi} being obvious for $n = 1$, and we first look for convenience at the cases when $n = 2$ or $n = 3$. We denote by $(e_1,e_2, \cdots ,e_n)$ a basis for $M$ consisting of elements of $\mathcal{S}$, and begin with a lemma that we shall use several times in the proof:

\begin{lemma}[\cite{MatPerf}, Lemma 6.1.6.] \label{Lemma 6.1.6.}
Under the hypotheses of Theorem \ref{Thm AnCombi}, we have: \par
1. If the components of two elements of $\mathcal{S}$ on two elements of $\mathcal{B}$ are nonzero, they are then equal or opposite. \par
2. There does not exist in $\mathcal{S}$ any $3$-tuple $(e_i - e_j, e_i - e_k, e_j + e_k)$.
\end{lemma}

\textit{Proof.} \; In the first case, after having changed if necessary the signs of $e_1$ and $e_2$, we would have the characteristic determinant
$$\operatorname{det}_{(e_i, e_j)}(e_i + e_j, e_i - e_j) = \begin{vmatrix}
1 & 1 \\ 1 & -1
\end{vmatrix} = -2;$$
similarly, we would have in the second case a third-order determinant
$$\operatorname{det}_{(e_i, e_j,e_k)}(e_i - e_j, e_i - e_k, e_j + e_k) = \begin{vmatrix}
 1 &  1 & 0 \\ 
-1 &  0 & 1\\
 0 & -1 & 1
\end{vmatrix} = 2.$$
\hfill $\square$ \\

If $n = 2$, a half-system of elements of $\mathcal{S}$ is $\{e_1 e_2, e_1 \pm e_2\}$. By Lemma \ref{Lemma 6.1.6.},
$e_1 + e_2$ and $e_1 - e_2$ cannot both belong to $\mathcal{S}$. Negating $e_2$ if necessary, we obtain the half-system $\{e_1, e_2, e_1 - e_2\}$ that we need. \par 

If $n = 3$, again by Lemma \ref{Lemma 6.1.6.}(1), there can be in $\mathcal{S}$ \textit{at most one} pair $\pm x$ with $x$ of the form $ \pm e_1 \pm e_2 \pm e_3$. By permutation and change of signs of $e_1, e_2, e_3$, we my assume that a half-system of elements of $\mathcal{S}$ consists of the five elements $e_1, e_2, e_3, e_4 = e_1 - e_2, e_5 = e_1 - e_3$ and a sixth one of the form $e_2 \pm e_3$ or $e_1 \pm e_2 \pm e_3$. Applying Lemma \ref{Lemma 6.1.5.}, we are left with the two possibilities $e_6 = e_2 - e_3$ or $e_6 = e_1 - e_2 - e_3$. In the first one we are done,
and the second case reduces to the first case by taking $e_2' = e_1 - e_2$ instead of $e_2$. \\

We now turn to the proof of Theorem \ref{Thm AnCombi} for $n \geqslant 4$. Let $M'$ be the submodule of $M$ generated by $e_2, e_3, \cdots, e_n$, and let $\mathcal{S}'$ be the set of nonzero elements $x \in M'$ such that $x$ or $\pm e_1 + x$ belongs to $\mathcal{S}$. Choose a half-system $\overline{\mathcal{S}}$ in $\mathcal{S}$ in such a way that the component on $e_1$ of any element of $\mathcal{S}$ be 0 or $+1$. One then obtains a map of $\overline{\mathcal{S}} \setminus \{ e_1 \}$ onto a half-system $\overline{\mathcal{S}'}$ of $\mathcal{S}'$ by suppressing the component on $e_1$ of the elements of $\mathcal{S}$. The inverse image of an element of $\overline{\mathcal{S}'}$ possesses one or two elements, this last possibility corresponding to twin systems $\{ x, e_1 + x \}$. To be able to use induction from $\mathcal{S}'$ to $\mathcal{S}$, we must prove the inequality $|\mathcal{S}'| \geqslant (n - 1)n$. To this end, we now prove the following lemma:

\begin{lemma}[\cite{MatPerf}, Lemma 6.1.7.] \label{Lemma 6.1.7.}
The number of twin systems $\{ x, e_1 + x \}$ in $\overline{\mathcal{S}}$ is at most $n-1$.
\end{lemma}

\textit{Proof. } Lemma \ref{Lemma 6.1.6.}(1) shows that two elements of $\mathcal{B}$ belonging to two twin systems may not have opposite nonzero components on one of the $e_i, i \geqslant 1$. Replacing $e_i$ by $-e_i$ when need be, we may assume that all their components are non-negative on the $e_i$. Suppose that there be at least $n$ twin systems, and let $x \in \mathcal{S}'$ such that $\{ x , e_1 + x \}$\footnote{In the textbook \cite{MatPerf}, it is writen as $(e_1, e_1 + x)$ there. This may be just a slip of the pen and I correct it in my paper here. } is a twin system and that $x$ has as many as possible components $x_i$ equal to 1. After permuting the $e_i$ if need be, we may assume that $x_2 = x_3 = 1$. There do not exist elements $y, z \in \mathcal{S}'$ belonging to distinct twin systems for which the components on $e_2$ and $e_3$ are respectively $(0,1)$ and $(1,0)$, for we would have
$$\operatorname{det}_{(e_1,e_2,e_3)}(x, e_1 + y, e_1 + z) = \begin{vmatrix}
 0 &  1 & 1 \\ 
 1 &  0 & 1\\
 1 &  1 & 0
\end{vmatrix} = 2.$$
We may thus assume that the system of components $(1,0)$ does not occur. Then replacing $e_2$ by $e_2 + e_3$\footnote{\label{indepOfTwin}Note: this element may not in $\mathcal{S}'$; fortunately, we do not need this condition, since we just need a basis of $M'$ to label the twin systems in some sense. } reduces the number of components equal to $1$ for $x$, and does not increase it for the members of the other twin systems. If the component of $x$ on $e_4$ is nonzero, we can in the same way consider the new basis obtained by replacing $e_2$ by $e_2 + e_4$, until each member of a twin system will possess a single nonzero component on the new basis. This contradicts the existence of at least $n$ twin systems, since $\mathcal{S}'$\footnote{Here shall be ${M}'$ rather than $\mathcal{S}'$, for the new basis we adopt may not and need not belong to $\mathcal{S}'$ as we explained in the former footnote. } contains only $n - 1$ elements of this type. \hfill $\square$ \\

\textit{End of the proof of \ref{Thm AnCombi}.} The induction hypothesis first shows that $\mathcal{S}'$ contains
exactly $(n-1)n$ pairs $\pm x$. Since one obtains $\mathcal{S}$ by adjoining to $\mathcal{S}'$ the
elements $\pm e_1$ and at most $n - 1$ pairs $\pm (e_1 + x)$ with $x \in \mathcal{S}'$, we see that $\mathcal{S}$ contains exactly $\frac{n(n+1)}{2}$ pairs $\pm x$, namely those of $\mathcal{S}'$, the pair $\pm e_1$, and
$n - 1$ pairs $\pm (e_1 + x)$ with $x \in \mathcal{S}'$.\footnote{\label{ErrorI} The first mistake in the proof. }

There remains to prove that $\mathcal{S}$ contains a basis for $M$ of the type we want. The induction hypothesis implies that up to a permutation of the elements of $\mathcal{S}'$, we may assume that $\mathcal{S}'$ contains a half-system of the form
$$ \overline{\mathcal{S}'} = \{\pm e_i \mid 2 \leqslant i \leqslant n \} \cup \{ \pm(e_i-e_j) \mid 2 \leqslant i < j \leqslant n \};$$
this implies that $\mathcal{S}$ contains the half-system $\overline{\mathcal{S}}$ which one obtains from $\overline{\mathcal{S}'}$ by
adjoining $e_1$ and of $n - 1$ elements of the form $e_1 \pm e_i$ or $e_1 \pm (e_i - e_j)$. 

Let $r$ be the number of elements of $\mathcal{S}$ of the first kind. We shall show that we can restrict ourselves by a suitable permutation of the elements of $\mathcal{S}'$ to
the case where $r = n - 1$.

We remark that the following transformations stabilize $\mathcal{S}'$ \footnote{\label{ErrorII} The second mistake in the proof. } and thus allow a partition of $\mathcal{S}'$ into two half-systems, one of which is of the desired form:
\begin{itemize}
\item The simultaneous change of signs of $e_2, \cdots, e_n$.

\item The permutations of $e_2, \cdots, e_n$.

\item The transformations $p_{i,j} (i, j \geqslant 2, i \neq j)$ which map $e_j$ onto $e_i - e_j$ and fix $e_k$ for $k \neq j$.
\end{itemize}

Since $p_{i,j}$ transforms $e_1 \pm (e_i - e_j)$ into $e_1 \pm e_j$, we may assume that $r \geqslant 1$. By a permutation of $e_2, \cdots, e_n$, we restrict ourselves to the case where the $r$ elements of the first kind are $e_1 \pm e_2, \cdots ,e_1 \pm e_{r+1}$ for a suitable choice of the signs. If there exists in $\mathcal{S}$ an element $e_1 \pm (e_i - e_j)$ with $i,j > r + 1$, applying $p_{i,j}$ allows us to increase $r$ to $r + 1$.

From now on, we suppose that such an element does not exist. Negating $e_2, \cdots, e_n$ if necessary, we may assume that there is a minus sign in front of $e_2$.
Lemma \ref{Lemma 6.1.6.} then shows that these $r$ elements are $e_1 - e_2, \cdots, e_1 - e_{r+1}$.

If $r < n - 1$, let $x = e_1 \pm (e_i - e_j)$ be another element of $\mathcal{S}$. We may assume that $i \leqslant r + 1$. Lemma \ref{Lemma 6.1.6.} shows that there must be a minus sign in front of $e_i$, hence that $x = e_1 - (e_i - e_j)$. Since there is a plus sign in front of $e_j$, we have $j > r + 1$.

Since the transformation $p_{i,j}$ preserve the elements $e_1 - e_k$ for $2 \leqslant k \leqslant r + 1$ and transforms $e_1 - (e_i - e_j)$ into $e_1 - e_j$, it can be used to increase again $r$ to $r + 1$. Iterating the process finally gives $r$ the value $n - 1$. \hfill $\square$ \\

\section{Some Mistakes and a Counter Example.} 

\quad \, There are several mistakes in the proof given above. One of them is easy to fix, while two of others are hard. 

\subsection{A Mistake Easy to Fix.}

\quad \, At the beginning of the proof of cases $n \geqslant 4$, the author constructed a set $\overline{\mathcal{S}'}$ by choosing a half-system of $\mathcal{S}$ and suppressing the component on $e_1$ of elements in $\overline{\mathcal{S}}\setminus\{e_1\}$. This may cause that $\overline{\mathcal{S}'}$ is not a half-system. In fact, if $\{\pm e_1, \pm(e_1 - e_2)\}$ is a twin system, one happens to choose $e_1 - e_2$ and $e_2$ when he determines $\overline{\mathcal{S}}$. Suppressing $\overline{\mathcal{S}}$, one may obtain that both $-e_2$ and $e_2$ belong to $\overline{\mathcal{S}'}$. \\ 

This phenomena is easy to overcome. We can exchange the order of constructing a half-system and suppressing the component on $e_1$. Namely, we can construct the set $\mathcal{S}' = \{ x \in M' \mid x \in M \; \text{or} \; x\pm e_1 \in M, x\neq 0 \}$ first, and next let $\overline{\mathcal{S}'}$ be an arbitrary half-system of $\mathcal{S}'$. This guarantees that $\mathcal{S}'$ will satisfy the induction conditions and $\overline{\mathcal{S}'}$ be a true half-system of $\mathcal{S}'$.

\subsection{Mistakes not trivial.}

\quad \, The first mistake in the original proof, which is marked by footnote \ref{ErrorI} in page \pageref{ErrorI}, is recovering the set $\mathcal{S}$ by adding $n-1$ pairs of $\pm (e_1 + x)$ to $\mathcal{S}'$. 

Although there are exactly $n-1$ twin systems in $\mathcal{S}$, which contribute $2(n-1)$ pairs of vectors to $\mathcal{S}$, the inverse images of other elements of $\mathcal{S}'$ may not be themselves. In fact, according to the definition of $\mathcal{S}'$, the inverse image of $x \in \mathcal{S}'$ could be $x$, $x+e_1$ or $x-e_1$. Therefore, we should recover $\mathcal{S}$ by not only adding $n-1$ twin system to $\mathcal{S}'$, but also modifying each vector $x \in \mathcal{S}'$ not in a twin system with adding $e_1$, $-e_1$ or 0 accordingly. In the following subsection \ref{secIll}, we can see that, for any integer $r$ less than $n$, one can easily construct a choice of initial basis $e_1, e_2, \cdots, e_n$ such that there are $(r-1)(n-r)$ pairs of elements of $\mathcal{S}'$ whose inverse images lie out of $M' = \spans_{\mathbb{Z}}\{e_2,\cdots, e_n\}$. \\

The second mistake, marked by footnote \ref{ErrorII} in page \pageref{ErrorII}, is that the transformations $p_{i,j}$ do not stabilize the forms of elements in $\mathcal{S}$. For example, it maps the vector $e_{\ell} - e_{j}\, (\ell \neq i)$, if it belongs to $\mathcal{S}$, to $e_{\ell} - (e_i - e_j)$, which has 3 non-zero components. As a consequence, we  eliminate a vector $e_1 - (e_i - e_j)$ with three non-zero components by $p_{i,j}$, but introduce a series new vectors $e_{\ell} - (e_i - e_{j})\, (\ell \neq i)$ of this form.

Hence, we cannot conclude at the end that we obtain a basis of $M$ with respect to which $\mathcal{S}$ has the desire form. We can only say that under the basis we obtain, all of $\pm e_i(1\leqslant i \leqslant n)$ and $\pm (e_1 - e_j) \, (2\leqslant j \leqslant n)$ belong to $\mathcal{S}$. Luckily, we can complete the proof based on the conclusion, although there is a more straightforward way to show this as presented in the proof of Lemma \ref{myLem1}.

\subsection{A Counter Example to Illustrate the First Fault.} \label{secIll}

\quad \, Here is a simple example which explains why we cannot obtain $\mathcal{S}$ just by adding $n-1$ new pairs of vectors to $\mathcal{S}'$. 

Let $f_1, \cdots, f_n$ be the $n$ independent vectors such that $\mathcal{S} = \{\pm f_i \mid 1 \leqslant i \leqslant n \} \cup \{ \pm(f_i - f_j) \mid 1 \leqslant i < j \leqslant n \}$. If we happen to choose the initial basis as
$$e_1 = f_1, \cdots, e_r = f_r, \quad e_{r+1} = f_{1} - f_{r+1}, \cdots, e_{n} = f_{1} - f_{n}.$$
Then, the vectors of $\mathcal{S}$ can be expressed under the basis we choose by
$$\begin{matrix}
\pm e_{1} & \pm (e_1 - e_2) & \cdots & \pm (e_1 - e_{r}) & \pm e_{r+1} & \cdots & \pm e_{n} \\
 & \pm e_{2} & \cdots & \pm(e_2 - e_r) & \pm (e_1 - e_2 -e_{r+1}) & \cdots & \pm(e_1 - e_2 - e_n)\\
 &          & \ddots & \vdots & \vdots & & \vdots \\
 &          &        & \pm e_{r} & \pm (e_1 -e_r - e_{r+1}) & \cdots & \pm(e_1 - e_r -e_n)\\
 &          &        &           & \pm (e_1 - e_{r+1}) & \cdots & \pm(e_{r+1} - e_n)\\
 &          &        &           &   & \ddots & \vdots \\
 &          &        &           &   &  & \pm(e_1 - e_n) \\
\end{matrix}.$$
In above table, the elements in main diagonal are $\pm f_i$, and the element in position $(i,j)$ is the difference $f_{i} - f_j$. \\

Hence, as was done in Section 2, the submodule $M' = \spans_{\mathbb{Z}} \{ e_2, e_3 \cdots, e_n \}$ contains 
$$ \mathcal{S}_1 = \{ \pm e_i \mid 2 \leqslant i \leqslant n\} \cup \{ \pm (e_i - e_j) \mid 2 \leqslant i < j \leqslant r \; \text{or} \; r+1 \leqslant i<j \leqslant n \}, $$
which are also contained in $\mathcal{S}$. To obtain the compressed set $\mathcal{S}'$, we need to combine $\mathcal{S}_1$ with the set
$$\mathcal{S}_2 = \{ \pm (e_i - e_j) \mid 2 \leqslant i \leqslant r, \; r+1 \leqslant j \leqslant n \},$$
the elements of which are not in $\mathcal{S}$. Then $\mathcal{S}' = \mathcal{S}_1 \cup \mathcal{S}_2$. \\

In this way, the $n-1$ twin systems are $\{ -e_i, e_1 - e_i \} _{2\leqslant i \leqslant n }$. But we cannot recover $\mathcal{S}$ by just adding these $n-1$ pairs of vectors to $\mathcal{S}'$. In fact, we should modify each element in $\mathcal{S}_2$ by adding or subtracting $e_1$ accordingly to pull it back to $\mathcal{S}$. We need to modify $|\mathcal{S}_2|/2 = (r-1)(n-r)$ pairs of vectors in total.

\section{The Corrected Proof.} \label{secCorPf}

\quad \, Since the main mistakes in the original proof are the construction of $\mathcal{S}$ from a basis of $\mathcal{S}'$ and the usage of transformations $p_{i,j}$, which have nothing to do with Lemma \ref{Lemma 6.1.6.}, Lemma \ref{Lemma 6.1.7.} and the assertion that $|\mathcal{S}| = (n+1)n$, we can give our corrected proof based on these propositions. \\

Actually, despite of the irregular form of the components of some element, we can conclude at the end of the original proof that:

\begin{lemma}\label{myLem1}
There exist $n$ independent elements $e_1, e_2, \cdots, e_n$ of $\mathcal{S}$ such that every $e_1 - e_i(2 \leqslant i \leqslant n)$ belongs to $\mathcal{S}$. 
\end{lemma}

\textit{Proof.} \; We can show this by using the transformations $p_{i,j}$ to increase $r$ repeatedly until $r = n-1$ as was done in the original proof. However, there is a more straightforward way to the conclusion.

By the induction hypothesis and Lemma \ref{Lemma 6.1.7.}, we have $|\mathcal{S}| = n(n+1)$,  which implies that there are exactly $n-1$ pairs of twin systems $\{x_1, e_1 + x_1\}, \cdots, \{x_{n-1}, e_1 + x_{n-1}\}$. Let $e_{i} = - x_{i-1}(2\leqslant i \leqslant n)$, and we will obtain that $e_1, e_2, \cdots, e_n$ belong to $\mathcal{S}$ and $e_1 - e_i(2 \leqslant i \leqslant n)$ also belong to $\mathcal{S}$. As for the independency of $e_2,\cdots,e_n$, we have seen in the proof of Lemma \ref{Lemma 6.1.7.}(also see footnote \ref{indepOfTwin} in Page \pageref{indepOfTwin}) that we can choose a suitable basis of $M'$ such that every $x_i$ has exactly one non-zero component under this basis. Hence, $x_1,\cdots, x_{n-1}$ form a basis of $M'$ and of course are independent. So are $e_2, \cdots, e_n$. \hfill $\square$ \\

Choosing the $n$ independent elements $e_1, e_2, \cdots, e_n$ obtained in the above lemma as a basis of $M$, we can characterize the remaining elements in $\mathcal{S}$ as following:

\begin{lemma} \label{myLem2}
Beside $\pm e_1, \pm e_2, \cdots, \pm e_n$ and $\pm (e_1 - e_i), 2 \leqslant i \leqslant n$ belonging to $\mathcal{S}$, other elements in $\mathcal{S}$ can only have two non-zero components with the form $\pm(e_i - e_j), 2 \leqslant i < j \leqslant n$ or three non-zero components with the form $\pm(e_1 -e_i - e_j), 2 \leqslant i < j \leqslant n$. Moreover, for each index pair $(i,j), 2 \leqslant i < j \leqslant n$, one and only one pair of $\pm(e_i - e_j)$ or $\pm(e_1 - e_i - e_j)$ belong to $\mathcal{S}$. 
\end{lemma}

\textit{Proof.} \; Let $y \in \mathcal{S}$ be a vector not listed in Lemma \ref{myLem1}. If the component of $y$ with respect to $e_1$ is 0, arbitrary two components of $y$ must have opposite signs. Otherwise one of $\pm y$ has component 1 on some $e_i$ and $e_j$ both, then we will have
$$\operatorname{det}_{(e_1, e_i, e_j)} ( y, e_1 - e_i, e_1 - e_j ) = \begin{vmatrix}
0 &  1 &  1 \\
1 & -1 &  0 \\
1 &  0 & -1
\end{vmatrix} = 2,$$
which is impossible. Therefore $y$ can only have two non-zero components, since there must exist two integers with the same sign among three or more non-zero integers. 

If the component of $y$ with respect to $e_1$ is non-zero, using $-y$ to substitute $y$ when necessary, we can assume the component is 1. In this case, the other non-zero components of $y$ have to be $-1$. Otherwise the characteristic determinant 
$$\operatorname{det}_{(e_1, e_i)} ( y, e_1 - e_i ) = \begin{vmatrix}
1 &  1  \\
1 & -1 
\end{vmatrix} = -2$$
would appear. Moreover, the number of non-zero components of $y$ could not be four or more, in which case we would have a forth-order characteristic determinant
$$\operatorname{det}_{(e_1, e_i, e_j, e_k)} ( y, e_1 - e_i, e_1 - e_j, e_1 - e_k ) = \begin{vmatrix}
 1 &  1 &  1 &  1 \\
-1 & -1 &  0 &  0 \\
-1 &  0 & -1 &  0 \\
-1 &  0 &  0 & -1
\end{vmatrix} = 2.$$
Hence, $y$ can only have the form $e_1 - e_i - e_j$.

Noticing that $e_i - e_j$ and $e_1 - e_i - e_j$ cannot belong to $\mathcal{S}$ both due to 
$$\operatorname{det}_{(e_i, e_j)} (e_i - e_j, e_1 - e_i - e_j) = \begin{vmatrix} 1 & -1  \\ -1 & -1 \end{vmatrix} = -2,$$ 
and that there are $\frac{(n-1)(n-2)}{2}$ index pairs $(i,j), 2 \leqslant i < j \leqslant n$, the cardinality of $\mathcal{S}$ is bound above by
$$\frac{1}{2} |\mathcal{S}| \leqslant n + (n-1) + \frac{(n-1)(n-2)}{2} = \frac{n(n+1)}{2}.$$
But $|\mathcal{S}| = n(n+1)$, which implies that one of $e_i - e_j$ or $e_1 - e_i - e_j$ belongs to $\mathcal{S}$ exactly for each index pair $(i,j), 2 \leqslant i < j \leqslant n$. \hfill $\square$ \\

With the help of above lemma, we can arrange the elements of $\mathcal{S}$ in a table shown following:
\begin{equation} \label{arrangeTable}
\begin{matrix}
\pm e_1 & \pm(e_1 - e_2) & \pm(e_1 - e_3) & \cdots & \pm(e_1 - e_{n-1}) & \pm(e_1 - e_n) \\
        & \pm e_2        & \pm y_{2,3} & \cdots & \pm y_{2,n-1} & \pm y_{2,n} \\
        &  & \pm e_{3} & \cdots & \pm y_{3,n-1} & \pm y_{3,n} \\
        &                &   & \ddots & \vdots & \vdots \\
        &                &   &   & \pm e_{n-1} & \pm y_{n-1,n} \\
        &                &   &   &   & \pm e_{n} \\

\end{matrix}
\end{equation}
Each $y_{i,j}(2 \leqslant i < j \leqslant n)$ has the form $e_i - e_j$ or $e_1 - e_i - e_j$. \\

\textit{The Proof of Theorem \ref{Thm AnCombi}. } Let us consider the last column in the table from $y_{2,n}$ to $y_{n-1,n}$. If $y_{j,n} = e_j - e_n$, we are done and consider the next $y_{j+1,n}$. If $y_{j,n} = e_1 - e_j - e_n$, we adopt a new basis in $\mathcal{S}$ as
$$ e_k' = e_k \, (k \neq j) ,\qquad e_j' = e_1 - e_j .$$
Then, we still have all of
$$e_1' - e_k' = e_1 - e_k \, (k \neq j), \qquad e_1' - e_j' = e_j$$
belong to $\mathcal{S}$. Therefore the new basis also satisfy the desired condition in Lemma \ref{myLem1}, and the elements in the first row of table \eqref{myLem2} do not change their forms in the new basis. Moreover, the element $y_{j,n} = e_1 - e_j - e_n$ has the form $e_j' - e_n'$.   

When we have done suitable substitutes of basis for all $j = 2, 3, \cdots, n-1$, we will eventually obtain a basis $e_1,\cdots,e_n$ with the elements in last column having the form $y_{j,n} = e_j - e_n$. 

Now, elements who may have three non-zero coordinates can only be located at $(i,j)$-place in the table, $2\leqslant i < j \leqslant n-1$. Fortunately, for each $y_{i,j}$, on the one hand, as was shown before, it can only have the form $e_i - e_j$ or $e_1-e_i-e_j$ since $\left\lbrace e_1 - e_i \mid 2 \leqslant i \leqslant n\right\rbrace  \subseteq \mathcal{S}$. On the other hand, $\left\lbrace e_j - e_n \mid 1 \leqslant j \leqslant n-1 \right\rbrace \subseteq\mathcal{S}$ forces $y_{i,j}$ can only be $e_i - e_j$ or $\pm(e_n - e_i - e_j)$. Therefore, $y_{i,j}$ must be $e_i - e_j$. In this way, we obtain $n$ independent elements $e_1, e_2, \cdots, e_n$ of $\mathcal{S}$ such that
$$\mathcal{S} = \left\lbrace \pm e_i \mid 1 \leqslant i \leqslant n \right\rbrace  \cup \left\lbrace  \pm(e_i-e_j) \mid 1 \leqslant i < j \leqslant n \right\rbrace .$$

\end{document}